\documentclass[12pt]{article}
\usepackage{amsmath,amssymb,amsfonts}
\usepackage{setspace}
\usepackage[margin=2cm]{geometry}
\usepackage{microtype}
\usepackage{booktabs}
\usepackage[final]{graphicx} 
\usepackage{color}
\usepackage{placeins}
\usepackage{rotating}
\usepackage{epstopdf}
\usepackage{subfig}
\usepackage{pdflscape}
\usepackage{enumerate}
\usepackage{float}
\usepackage{longtable}
\usepackage{authblk}
\usepackage{caption}
\usepackage{textcomp}
\usepackage[utf8]{inputenc}
\usepackage{algorithm}
\usepackage{algpseudocode}
\usepackage[colorlinks=true, linkcolor=blue, citecolor=blue]{hyperref}
\usepackage{natbib}
\usepackage{alltt}
\usepackage[labelsep=period]{caption}
\usepackage{tikz}
\usetikzlibrary{shapes.geometric,arrows.meta}
\usetikzlibrary{positioning,calc}
\usepackage[symbol]{footmisc}
\usepackage{multicol}
\usepackage{rotating}

\title{On Improvement of Control Chart using Repetitive Sampling for Monitoring Process Mean}

\author[1,3*]{Fahad Rafique}
\author[3]{Saadia Masood}
\author[2]{Shabbir Ahmad}
\author[1,3]{Sadaf Amin}

\affil[1]{School of Mathematical Sciences, Capital Normal University Beijing, China}
\affil[2]{SDepartment of Mathematics, Comsats university wah Campus, Pakistan}
\affil[3]{Department of Statistics, PMAS-Arid Agriculture University, Rawalpindi,  Pakistan}
\date{}
\begin{document}
	\maketitle
\begin{abstract}
	In the practical industry, the most commonly used application of statistical analysis for monitoring the process mean is the control chart. Control charts are generated based on the presumption that we have a sample from a stable process. The control chart then provides a graphical display to test this presumption. In the existing estimator \textcolor{red}{$Mr$}, researchers use a technique involving repetitive sampling along with an auxiliary variable for detecting and monitoring the statistical process mean. The existing control chart, namely \textcolor{red}{$Mr$}, is based on the regression estimator of the mean using a single auxiliary variable $X$. We propose the \textcolor{red}{$Mrep$} chart using a ratio-product exponential type estimator, and the \textcolor{red}{$Mrwp$} chart with a more efficient difference-cum-exponential type estimator used in quality control for improving the process mean in terms of $ARL$. Then we compare the proposed charts \textcolor{red}{$Mrep$} and \textcolor{red}{$Mrwp$} with the existing \textcolor{red}{$Mr$} chart in terms of $ARL$. Using $ARL$ as a performance measure, better results of the proposed charts are observed for detecting shifts in the mean level of the characteristic of interest. Moreover, Monte Carlo simulation in terms of repetitive sampling is used for quality control charting and statistical process control for the betterment of the process mean.
\end{abstract}


\maketitle

\section{Introduction}
\label{sec:introduction}
Quality control (QC) is a structure of everyday technical activities, to measure and control the quality of the inventory as it is far being advanced. The quality control tool is designed to deliver ordinary and consistent exams to ensure statistics integrity, correctness, and totality. Process monitoring has a very ancient history. Shewhart \cite{Shewhart1924} presented the procedure for monitoring of process mean. Moreover, in the field of process monitoring Statistical Process Control became a vital tool. In statistical process control methods are observed for their constancy in terms of different parameters by using statistical methods like histograms, check sheets, pareto diagrams and control charts.

For a better monitoring and evolution of a quality characteristic of interest, \cite{Montgomery2009} defines SPC as a controlling gathering of problem-solving implements advantageous in attaining improving capability and process stability over then done with the decline of variation. \cite{Isaki1983} studied that variance estimators under many sample designs were proposed and compared once auxiliary data was offered. \cite{LuReynolds1999} suggested that in statistical process control it was generally assumed that the observations picked from the process are independent, but in routine the observations in many cases are actually auto correlated. \cite{Petcharat2023} suggested that HEWMARS-CV (hybrid exponentially weighted moving average control chart) control chart used the general repetitive sampling method to generate two parallel pairs of control limits to boost the performance of control charts. \cite{Aslam2019} discuss the application of the exponentially weighted moving average (EWMA) control chart for the monitoring blood glucose in type-II diabetes patients. A \cite{Saghir2019} proposed that exponentially weighted moving average (EWMA) control charts using repetitive sampling technique as well as an auxiliary variable for very efficient detection of small to moderate shifts in location type control charts. In Generalized family of EWMA charts \cite{Lazariv2015} state that three different repetitive sorts of sampling plans. \cite{Aslam2019} suggest that for some more literature, in these same directions, one may view \cite{ChengChen2010}, \cite{Riaz2008}, \cite{SantiagoSmith2013}, \cite{Ahmad2014}, \citealp{Shabbir2014} and \cite{Yaqub2016}.

\cite{Riaz2008} suggested that a control chart, which is basically the type of Shewhart-type control chart, was established for detecting the shift of process mean level, and we know the basic function of this chart is detecting the large shift more quickly. Particularly the chart depends on the auxiliary variable X, which is relayed on the regression estimator of mean. Considering the assumption of bivariate normality, as he assumed that the standard normal bivariate random variable for the layout of the chart was evolved for section I exceptional management. For similar motive proposed charts are contrasted with the existing control chart. For overall performance of degree he was using the power curve, and using that power curve he made improvement in the performances of the chart, and he also noticed that for characteristic of interest this method was given change for detecting the shifts in mean level.

\section{METHODOLOGY}
They use the information given in the study variable or variable of interest with addition of auxiliary variables to improve the precision alongside with which parameter are estimated. Moreover information provided is depending on auxiliary variable and study variable have a strong. \cite{Riaz2008} introducing the concept of improving the process monitoring in control chart with the impression of using the information based on an auxiliary variable with its characteristics. For the improvement of monitoring the variance, he was using the information of quality characteristics of X as auxiliary variable with study Y. for the betterment in the improvement, in monitoring of variance of Y he was take benefit from the correlation between  and for variance of Y he proposed a usual regression-type estimator and acclaimed its usage in control charts. He was claiming the supremacy of his $M_r$ chart to the existing $S^2$ control chart, with condition on $\lvert\rho_{xy}\rvert$ [wherever $\rho_{xy}$ is the linear correlation between X and Y].

In quality control literature there are different categorization or classifications of control charts, for example   according to the type of data, sample size, type of control chart, and shape of control chart. \cite{Farnum1994} categorized two elementary kinds of control: deviation control and threshold control. Deviation control are concerned with detecting the small shifts in the process parameters whether for detection of large shift the Threshold control are concerned. Non-Shewhart control charts, (e.g. EWMA as well as  CUSUM charts) are regarded as deviation control charts while The Shewhart type control charts are regarded as threshold control charts. Following \cite{Riaz2008} he improved monitoring of the process mean level of a quality characteristic of interest Y, using information about a single auxiliary characteristic X. he proposed a Shewhart-type process mean control chart, namely   chart (is type of a threshold control charts) with considering the assumption of bivariate normality of $(Y,X)$. The main emphasis of the suggested chart is based on phase I level of quality control. Using a particular auxiliary variable X with the regression estimator for mean of Y (the study variable), that are clearly defined for a bivariate random sample, for example $(y_1,x_1), (y_2,x_2), \cdots (y_n,x_n)$ of size n as:
\begin{equation}
	M_r=\bar{y}+b(\mu_x-\bar{x}),
\end{equation}
where $\bar{y}, \bar{x}$ are the sample mean Y and X respectively, $\mu_x$ is the population mean of $X$ (consider to be known) and $b$ is well-defined as:
\begin{equation}
	b=r_{xy}(\frac{s_y}{s_x}),
\end{equation}
wherever $r_{xy}$ is the sample correlation coefficient among study variable Y and auxiliary variable X, moreover $s_y$ is the sample standard deviation of Y and $s_x$ is the sample standard deviation of X.

The basic three parameters of any Shewhart-type control chart are the lower control limit $(LCL)$, central line $(CL)$ and upper control limit $(UCL)$. There was basically two approaches to define these parameters, that is to say the 3-sigma limits approach and the probability limits approach. If the distributional performance of the related estimator is approximately symmetric, then we say that 3-sigma limits approach is a worthy alternative. In case of asymmetric distributional performance of the related estimator, the probability limits approach is preferred. The parameters of the existing $M_r$ chart using together the approaches are expressed as in equation (1).

The value  corresponding to control Limit of the existing $M_r$ chart, just like $\bar{y}$ for $R$ chart providing in \cite{Alwan2000} for supposing the probability of producing an type-I error to be fewer than a specifying value say $\alpha$, the control limits or we can say that actually true probability limits.

For the existing $M_r$ chart were define below.
\begin{equation}
	LCL=M_{r1} ~~\text{with}~~ p_n(M_r=M_{r1}) \le \alpha_1
\end{equation}
\begin{equation}
	UCL=M_{r\mu} ~~\text{with}~~ p_n(M_r=M_{r\mu}) \ge 1- \alpha_\mu
\end{equation}
where $\alpha=\alpha_1+\mu$ and $p_n$ representing the cumulative distribution function for a present value of n.
\begin{equation}
	LCL=M_{r1}=\bar{M}_r+\frac{c_l\bar{\sigma}_y}{\sqrt{n}} ~~\text{with}~~ p_n(C=C_1) \le \alpha_1,
\end{equation}
\begin{equation}
	UCL=M_{r\mu}=\bar{M}_r+\frac{c_\mu\bar{\sigma}_y}{\sqrt{n}} ~~\text{with}~~ p_n(C=C_\mu) \ge 1- \alpha_1.
\end{equation}
Hence, the quintile points of the distribution of C, the average of sample $M_rs$   this will $M_r$ and $\bar{\sigma}_y$ (which is defines as controlled estimate of the process standard deviation $\sigma_y$) allow setting the true probability limits for the existing $M_r$ chart.

If normal approximation to the distribution of $C$ is used, as we know C is constant then the usual 3-sigma control limits with the parameters of  $M_r$ Chart are given as:

\begin{align}
	&\
	\left.\begin{aligned}
		LCL & = \bar{M}_r -3\sigma_{Mr} \\
		CL & = \bar{M}_r \\
		UCL & = \bar{M}_r +3\sigma_{Mr}
	\end{aligned}\right\}       
\end{align} 
Using (4) in (7) gives the following result are:
\begin{align}
	&\
	\left.\begin{aligned}
		LCL & = \bar{M}_r -3k_2\hat{\sigma}_y / \sqrt{n} \\
		CL & = \bar{M}_r \\
		UCL & = \bar{M}_r +3k_2\hat{\sigma}_y / \sqrt{n}
	\end{aligned}\right\}       
\end{align}
where the values of $k_2$ are known as standard error. Consider that the how closely the normal approximation is to the true distribution of C, define the strength of these 3-sigma limit-based parameters of the existing. 

When we decided the structure of the control charts for the probability limit or either 3-sigma limit approach for the given level of significance, then we draw the sample statistics in favor of the sample statistic with time order of samples now then we conclude that if all the sample $M_r$ is lie within the control limits, that is clear evidence that no shift is used in the process mean level and $M_r$ is stable at the time. While if the process is unstable then is some assignable cause is present in the process and creating a shift at in the process at the level of process mean level. 

To address small and moderate shifts, using the developed structure of $M_r$ chart, Aldrich and Nelson \cite{5} given basic structure for  the run rules, to  produced small and moderate shift in the process with the help of developed structure of $M_r$ chart  the risk of false alarms is increased.

\section{RESULTS AND DISCUSSION}
\subsection{Proposed Technique 1}
To improve the process mean of control chart we followed a ratio product exponential type estimator $\hat{\bar{Y}}_{S,RP}$ of $\bar{Y}$ is suggested by \cite{SinghChauhanSawan2008} as:
\begin{equation}
	\hat{\bar{Y}}_{S,RP}=\bar{y}\bigg(\alpha ~ exp(\frac{\bar{X}-\bar{x}}{\bar{X}+\bar{x}})+ (1-\alpha)~ exp(\frac{\bar{x}-\bar{X}}{\bar{x}+\bar{X}})\bigg),
\end{equation}
where $\alpha$ is an arbitrary constant and optimum value of $\alpha$ as given, 
\begin{equation}
	\alpha_{opt}=\frac{1}{2}+\frac{\rho_{xy}C_y}{C_x}.
\end{equation}
Mean square error up to 1st order approximation which is accurately equal to $Var(\hat{\bar{Y}}_{Reg})$ is given as:
\begin{equation}
	MSE_{min}=\hat{\bar{Y}}_{S,RP}\approx\lambda \bar{Y}^2(1-\rho_{xy}^2)C_y^2=Var(\hat{\bar{Y}}_{Reg}).
\end{equation}

\subsection{Proposed Technique 2}
To improve the control chart using repetitive sampling for monitoring process mean we followed a more improve ratio type estimator including weights under the perspective of estimator proposed by \cite{Shabbir2014} as:
\begin{equation}
	\hat{\bar{Y}}^*_{P}=\bigg(\hat{\bar{Y}}_{S,RP}+w_1(\bar{X}-\bar{x})+w_2\bar{y}\bigg) exp(\frac{\bar{X}-\bar{x}}{\bar{X}+\bar{x}}),
\end{equation}
as $\hat{\bar{Y}}_{S,RP}$ is given as:
\begin{equation}
	\hat{\bar{Y}}_{S,RP}= \alpha ~ exp(\frac{\bar{X}-\bar{x}}{\bar{X}+\bar{x}})+ (1-\alpha)~ exp(\frac{\bar{x}-\bar{X}}{\bar{x}+\bar{X}},
\end{equation}
the optimum values of $w_1$ and $w_2$ are given by
\scriptsize
\begin{equation}
	w_{1opt}=\frac{\bar{Y}\bigg(-4\rho_{xy}C_y+C_x(2-\lambda C_x^2+\lambda \rho_{xy}C_y C_x+2\lambda(-1+\rho_{xy}^2))C_y^2\bigg)}{4\bar{X} C_x\bigg(-1+\lambda (-1+\rho_{xy}^2)C_y^2\bigg)}
\end{equation}
\normalsize
\begin{equation}
	w_{2opt}=\frac{\lambda (C_x^2-4(-1+\rho_{xy}^2)C_y^2}{4\bar{X} C_x\bigg(-1+\lambda (-1+\rho_{xy}^2)C_y^2\bigg)}.
\end{equation}

\begin{equation}
	MSE_{min}(\hat{\bar{Y}}^*_P) \approx \frac{\lambda \bar{Y}^2\bigg(\lambda C_x^4-8((-1+\rho_{xy}^2)(-2+\lambda C_x^2)C_y^2)\bigg)}{16\bigg(-1+\lambda (-1+\rho_{xy}^2)C_y^2\bigg)}.
\end{equation}
After simplified (16) it is written as:
\begin{equation}
	MSE_{min}(\hat{\bar{Y}}^*_P) \approx MSE_{min}(\hat{\bar{Y}}_{Reg}) - (T_1+T_2).
\end{equation}
The optimum value of $T_1$ and $T_2$ are given below,
\begin{equation*}
	T_1=\frac{\lambda^2 \hat{\bar{Y}} \bigg(C_x^2+8(1-\rho_{xy}^2)C_y^2\bigg)^2}{64\bigg(1+\lambda (1-\rho_{xy}^2)C_y^2\bigg)},
\end{equation*}

\begin{equation*}
	T_2=\frac{\lambda^2 \hat{\bar{Y}} C_x^2 \bigg(3C_x^2+16(1-\rho_{xy}^2)C_y^2\bigg)^2}{64\bigg(1+\lambda (1-\rho_{xy}^2)C_y^2\bigg)}.
\end{equation*}

For monitoring quality characteristics of industrial or service systems, control charts are best graphical tools. Control charts had been implemented in different industrial organisms such as short-run production and long-run production ever since 1930. It has no more that too much fewer research work is done to implementing control charts in job-shop industrialized and batch production. In this thesis, computerized numerical simulations are used to find appropriate control charts for improvement of process mean. Monte Carlo simulation technique is used to find the value of ARL against different value of shift. The main task of this thesis is to study control charts capabilities to early detect the mean shifts for process mean production.

For over convenience we use abbreviation $T_0$, $T_1$, $T_2$ and $T_3$ in $R$ console statistical software for our existing and propose estimator.
\begin{equation*}
	T_0=\bar{y}
\end{equation*}
\begin{equation*}
	T_1=M_r=\bar{y}+b(\mu_x-\bar{x})
\end{equation*}
\begin{equation*}
	T_2=M_{rep}=\hat{\bar{Y}}_{S,RP}=\bar{y}\bigg(\alpha ~ exp(\frac{\bar{X}-\bar{x}}{\bar{X}+\bar{x}})+ (1-\alpha)~ exp(\frac{\bar{x}-\bar{X}}{\bar{x}+\bar{X}})\bigg),
\end{equation*}
\begin{equation*}
	T_3=M_{rwp}=\hat{\bar{Y}}^*_{P}=\bigg(\hat{\bar{Y}}_{S,RP}+w_1(\bar{X}-\bar{x})+w_2\bar{y}\bigg) exp(\frac{\bar{X}-\bar{x}}{\bar{X}+\bar{x}}).
\end{equation*}

In Figures 2 to 6 (a) to (aa) we plot ARL curve of $T_0$, $T_1$, $T_2$ and $T_3$ estimators based on bivariate normal distribution when $ARL = 200, 371, 500$ $n=5, 10, 15$ and $l_{xy}=0.30,0.60,0.90$. When we plot ARL against different values of shift, we observed that our proposed estimators $T_2$ and $T_3$ are more efficient than other existing charts $T_0$ and $T_1$.  

In  order to illustrate the practical application of the statistical control structures under our research  and our study highlights their importance in an highly efficient detection of changes in process parameters, we provide here theoritical as well as descriptive examples to compare the performance and efficiency  of the structures $T_0$ vs. $T_3$ (which is based on Y and Y and  X respectively) The variables that we used in  our study  are Y and X may refer to different real life variables such as (1)  in electrical  power generation and production from coal, let we may consider that amount of total electrical power generated as Y(megawatts) , normal  air temperature, and humidity in air as X; (2) to evaluate the computer pipelining, moreover  we may consider the speed up as (Y) and number of stages as (X) ; (3) also in high quality assurance of pharmaceutical products use in medical field, we may consider that  pharmaceutical products as (Y) units and temperature as (X); more example  (4) in the processes monitoring of production of electrical steel wire use in industry  we may consider the tensile strength as Y (psi) and the outside wire diameter as X (mm); (5) in hospital for  patient heart doctor to  monitoring Electrocardiography (ECG)  of heart, the  progression of a  electric wave produced on ECG machine can be considered as Y (per minute) and  the blood pressure as X , etc. For the study   purposes, we have generated data sets containing 50 subgroups using sample size n = 5, 10 and 15 from the bivariate normal distribution.

In-control situation we generate the first 30 observations , i.e., $\delta=0$, whereas  an out-of-control situation  the remaining 20 observations are generated with $\delta=1$ for both type of data sets. Based on these type  data sets, the computed  as wellas calculated  values of $T_0$ vs. $T_3$ are given in Tables 1-2. The control charting structures based on $T_0$ vs. $T_3$ are worked out for these data sets and their corresponding graphical displays are shown in Figure 1 (i-v). The symbols used in these figures are described as: $LCL_i$ and $UCL_i$ which represent the lower and upper control limitsof control chart respectively, for the structures based on $T_0$ vs. $T_3$ the solid lines refer to the charting statistic of $T_0$ and $T_1$, the dotted lines refer to the charting statistics of $T_2$ and $T_3$, while the dashed line refers to central limit computed from the sample data.

It is evident from Figure 1(i) that $T_0$ and $T_1$-based control structure detects 5 and 16 out-of-control signals at sub group sample number 30, 43, 47,48, 49 and 11, 31, 32, 33, 35, 36, 37, 39, 41, 42, 43, 44, 46, 47, 48, 49 while $T_2$ and $T_3$-based control structure detects 17 and 20 out-of-control signals at sample number 31, 32, 34, 35, 36, 37, 39, 40, 41, 42, 44, 45, 46, 47, 48, 49 and 31, 32, 33, 34, 35, 36, 37, 38, 39, 40, 41, 42, 43, 44, 45, 46, 47, 48, 49, 20. It means that $T_2$ and $T_3$-based control structure has given more out of control signals as compared to T0 and $T_1$-based control chart. Moreover, as the shift occurred after sample number $30, T_2$ and $T_3$ chart reacted quickly (first out-of-control signal is given after 30 samples, i.e. at sample number 31) as compared with the usual $T_0$ and $T_1$ chart (first out-of-control signal is given after 30 samples, i.e. at sample number 43) (Figure 1(i)). 

A similar higher detection ability of $T_2$ and $T_3$-based structure relative to $T_1$-based structure may be seen in Figure 1. One can also observe from Figure 1(ii-iii) that $T_3$ chart has more quick detection ability as compared to $T_2$ chart, which is in accordance with the findings of WECO. Similarly, we can observe the superiority of $T_2$ and $T_3$ charts from Figure 1(iv-v) which are based on the data-set using n = 5 and 10. It is observed that the detection ability of these charts is improved with the increment of sample size n as expected, for instance $T_2$ and $T_3$ charts detected 18 and 19 out-of-control signals, respectively, for n = 10 (Figure 1(v)). A similar high detection ability of the $T_3$ chart may be observed in the same figures.

\begin{figure*}[h]
	\centering
	\begin{tabular}{c c}
		\includegraphics[width=8cm, height=6cm]{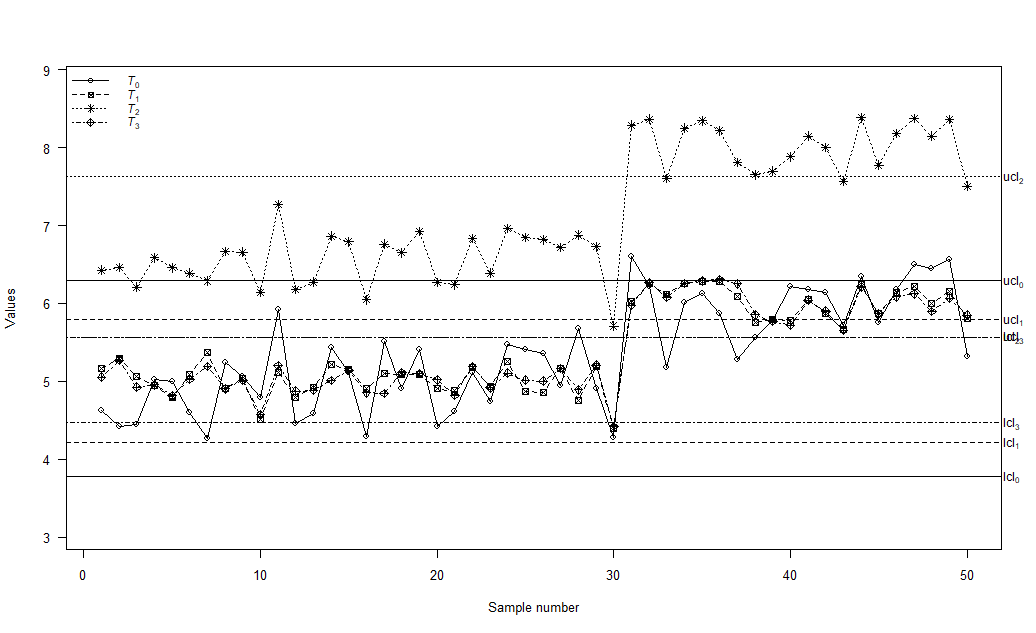}\hspace{0.1cm}&
		\includegraphics[width=8cm, height=6cm]{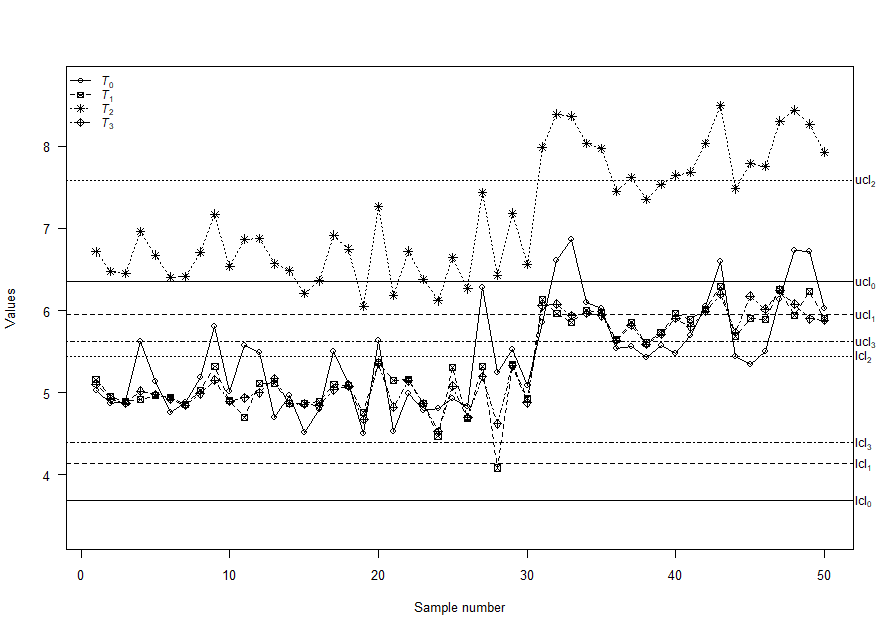}\hspace{0.1cm}\\
		(i) & (ii) \\
		\includegraphics[width=8cm, height=6cm]{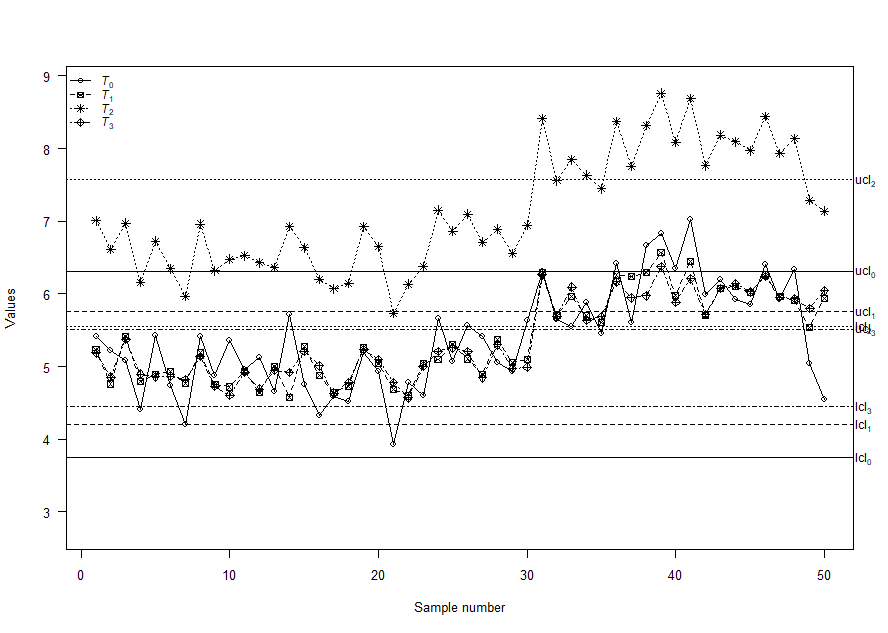}\hspace{0.1cm}&
		\includegraphics[width=8cm, height=6cm]{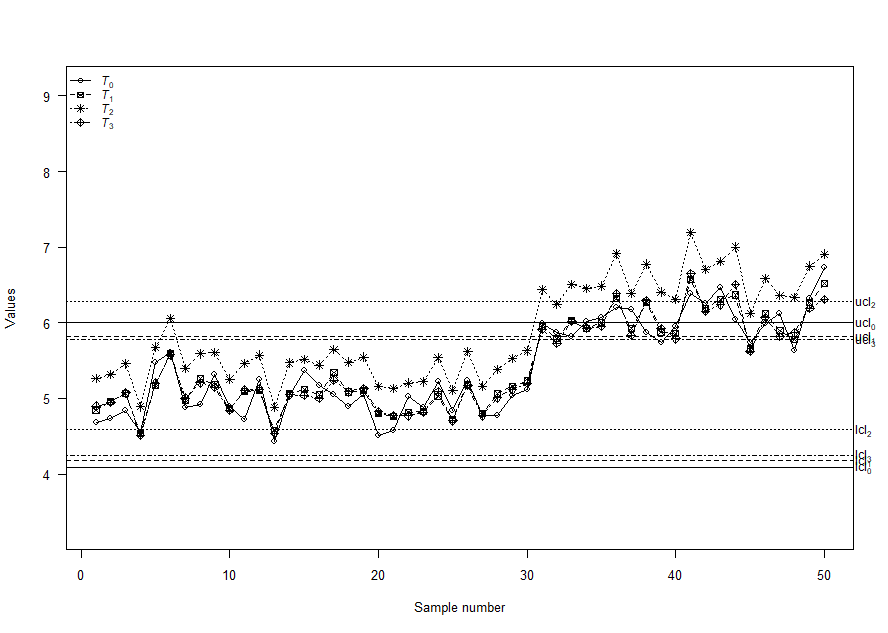}\hspace{0.1cm}\\
		(iii) & (iv) \\
	\end{tabular}
	\includegraphics[width=8cm, height=6cm]{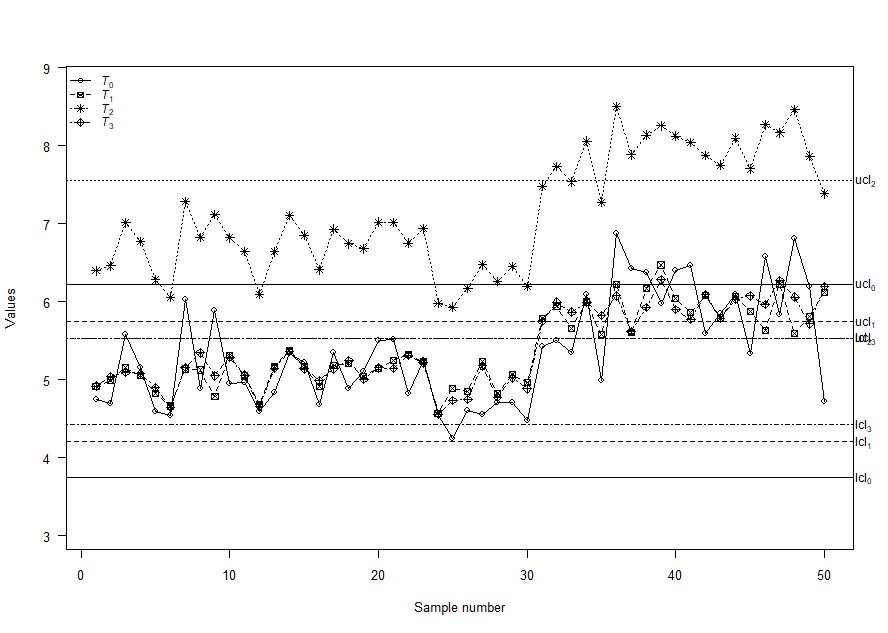}\hspace{0.1cm}\\
	(v)  \\
	\caption{$T_0-T_3$ based charts: (i) n=5, ARL=500, $\rho_{xy}=0.90$; (ii) n=5, ARL=371, $\rho_{xy}=0.90$; (iii) n=5, ARL=200, $\rho_{xy}=0.90$; (iv) n=10, ARL=371, $\rho_{xy}=0.60$; (v) n=15, ARL=500, $\rho_{xy}=0.90$.}
	\label{fig1}
\end{figure*} 

\section{CONCLUSION}
\label{sec:guidelines}
In every industrial manufacturing process there are two category of variation present namely which is naturally inherent in the process and other one is called the unnatural variation which is due to some special cause or defective manufacturing process. Three kind of unnatural variation are small, moderate and large shifts. We are using the shewhart's control charts in order  to detect large shift in the whole process as we know Shewhart control is used to detect large magnitude shift. This study is concerned with the improvement of control charts using auxiliary variables for monitoring process mean. $ARL$, $EQL$, $RARL$, and $PCI$ are computed using simulation technique with different value of shift. An extensive simulation study is conducted for each existing and purpose estimator to compare the efficiency of each. The purpose of simulation study is to analyze, on which estimator the ARL value is less or high and perform better in their respective in-control as well as out of control $ARL$. The simulation study for each estimator suggests that $T_0$ and $T_1$ estimator is less efficient than proposed estimators in term of $ARL$.

The proposed control charts tell us in detail about deterioration and improvement simultaneously in this study and consider that our purpose estimators give better improvement of control charts to ensure an efficient and proficient monitoring of process mean. The performance assessments are carried out in terms of four evaluation parameter $ARL$, $EQL$, $RARL$, and $PCI$. In general the rule of thumb is that we have a smaller EQL value for an efficient charting structure for process monitoring. In all cases When we calculated EQL value for $T_0-T_3$ We see that $T_3$ has smallest value as compare to other estimator. So $T_3$ is efficient charting structure for process mean. In same manner for RARL and PCI our proposed estimator give better performance.

Moreover we know that Control charts are very effective as well as useful tool to control and monitor the performance of any industrial, manufacturing and non-manufacturing process in term of quality control. Control charts that are based on auxiliary information are very effective and attractive alternatives to boost up the performance ability of their assembling arrangements and design structures. In this research study we presents different efficient control structures in the form of $T_0, T_1, T_2$ and $T_3$ charts (using different ways of using auxiliary information, including regression, ratio-product exponential type regression and difference-cum-exponential type estimator for improvement of process mean. It is detected that the control charts based on ratio-product exponential type regression ($T_2$) and difference-cum-exponential type estimator ($T_3$) charting statistics exhibited superior performance as compare to other charts in general. The T0 and T1 charts conveyed relatively lower performance as compare to other charting sutures $T_2$ and $T_3$ in general. For the prospective research analysis, the supposed design structures have the ability to execute or perform better than many other existing variability charts like $M_r, S$ and $R$ charts. The execution and implementation of this proposed design structures with the practical data sets is very simple, quick and attractive in detecting the shift and timely receiving the out-of-control signals.

An illustrated example is also conducted for check the performance of proposed charts. As $T_2$ and $T_3$ detect the shift more quickly and more signal point are detect in shift as compare to $T_0$ and $T_1$. From figure (1). we are observed that in out of control situation detect of point for $T_2$ and $T_3$  give efficient result as compare to existing estimators. The illustrative examples that is use in our study have also supported the foremost and major relevant findings of the research and shown that the quickly detection ability of shift and charting structures is improved/better with the increment of sample size.

The efficacy of proposed estimators for monitoring the process mean is obvious from the graphical display as by the performance measures of $ARL$, $EQL$, $RARL$, and $PCI$.

Further we recommended that simulation studies conducted for individually of the purposed estimator with shift give precise result as compare to without shift. This thesis is explored for so many other different real life distributions like Bayesian, and Multivariate distribution. The scope of the thesis may also be extended by: incorporating the information on more auxiliary variables, implementing runs rules, investigating varying sampling strategies, and generous recollection property to the proposed structures in the form of (EWMA and CUSUM) frames.

\begin{figure*}[h]
	\centering
	\begin{tabular}{c c}
		\includegraphics[width=8cm, height=6cm]{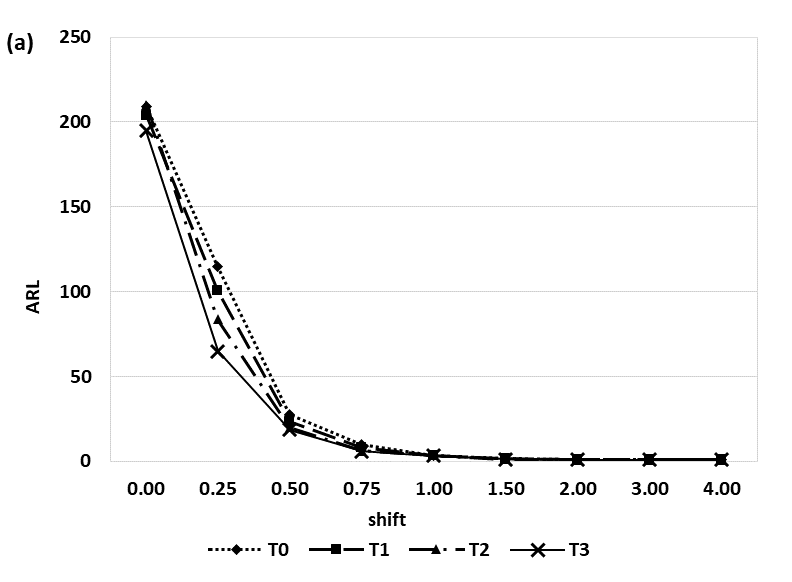}\hspace{0.1cm}&
		\includegraphics[width=8cm, height=6cm]{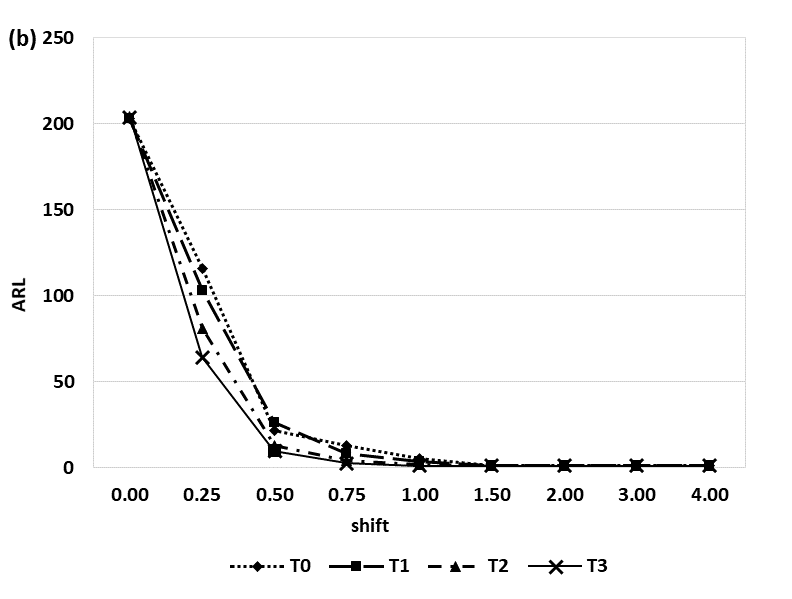}\hspace{0.1cm}\\
		
		\includegraphics[width=8cm, height=6cm]{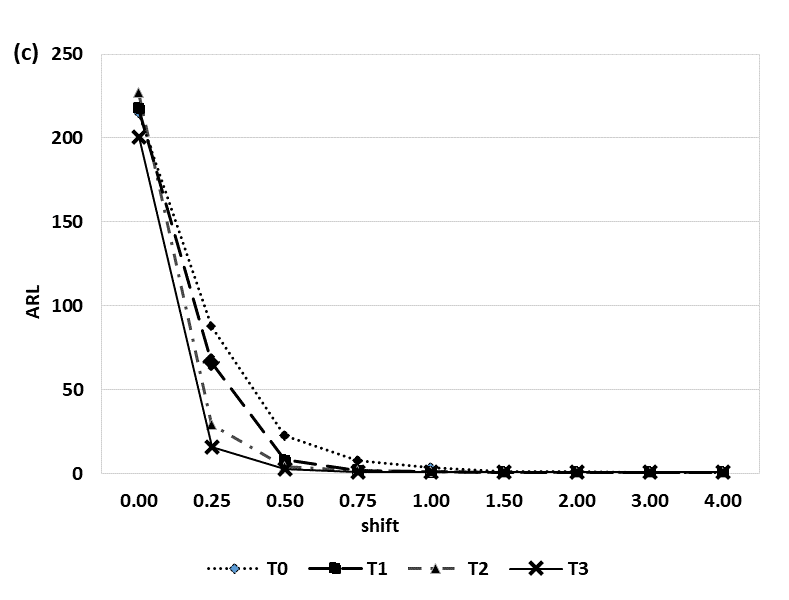}\hspace{0.1cm}&
		\includegraphics[width=8cm, height=6cm]{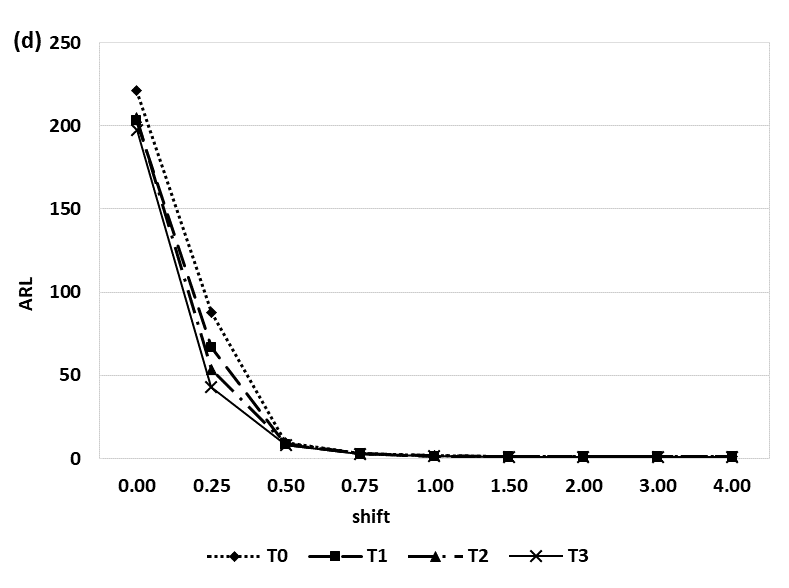}\hspace{0.1cm}\\
		
		\includegraphics[width=8cm, height=6cm]{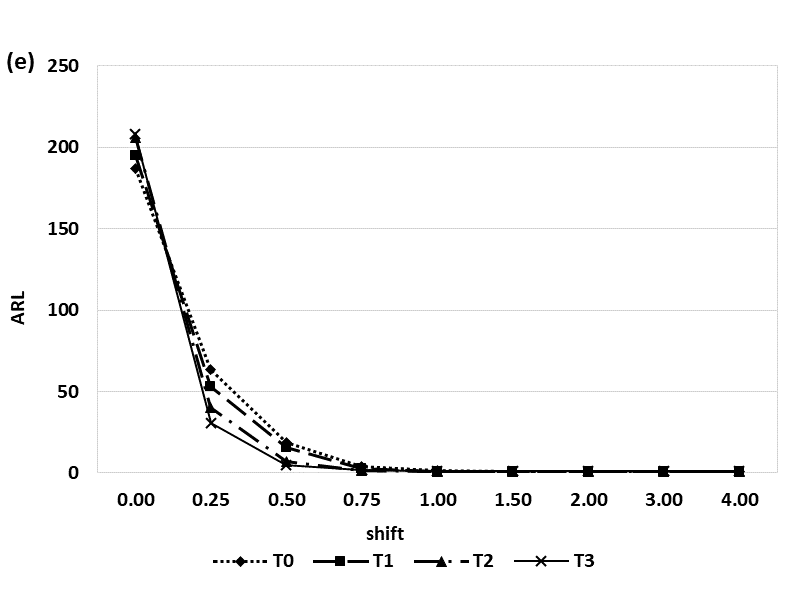}\hspace{0.1cm}&
		\includegraphics[width=8cm, height=6cm]{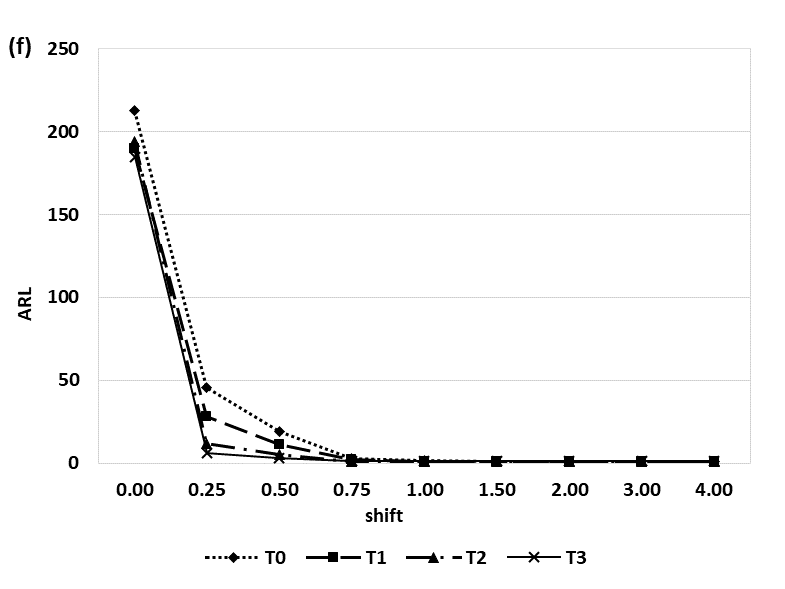}\hspace{0.1cm}\\
		
	\end{tabular}
	\caption{ARL curve of $T_0-T_3$ at: (a) ARL=200, n=5, $l_{xy}=0.30$; (b) ARL=200, n=5, $l_{xy}=0.60$; (c) ARL=200, n=10, $l_{xy}=0.90$; (d) ARL=200, n=10, $l_{xy}=0.30$; (e) ARL=200, n=10, $l_{xy}=0.60$; (f) ARL=200, n=10, $l_{xy}=0.90$.}
	\label{fig2}
\end{figure*} 

\begin{figure*}[h]
	\centering
	\begin{tabular}{c c}
		\includegraphics[width=8cm, height=6cm]{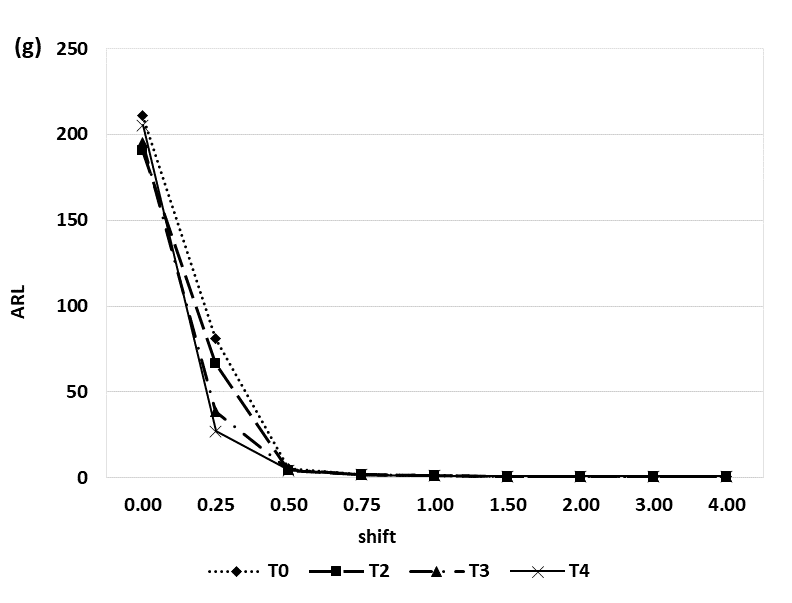}\hspace{0.1cm}&
		\includegraphics[width=8cm, height=6cm]{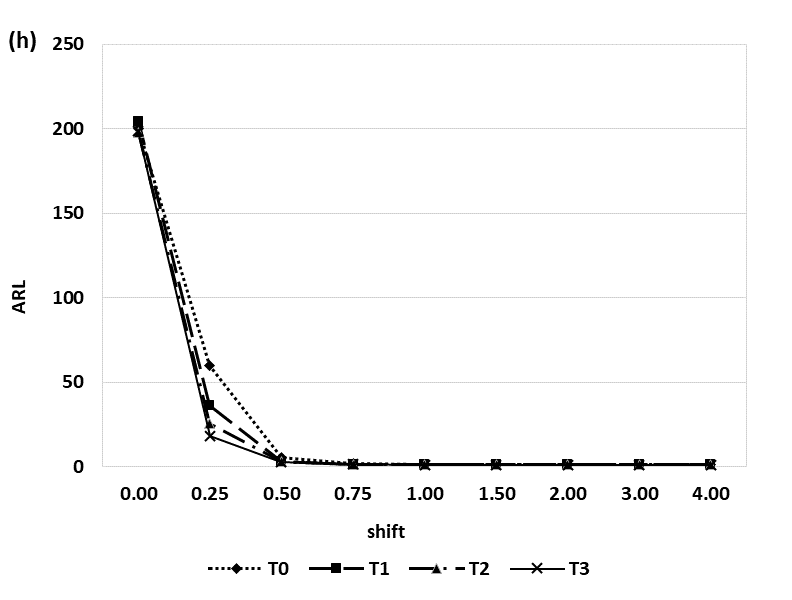}\hspace{0.1cm}\\
		
		\includegraphics[width=8cm, height=6cm]{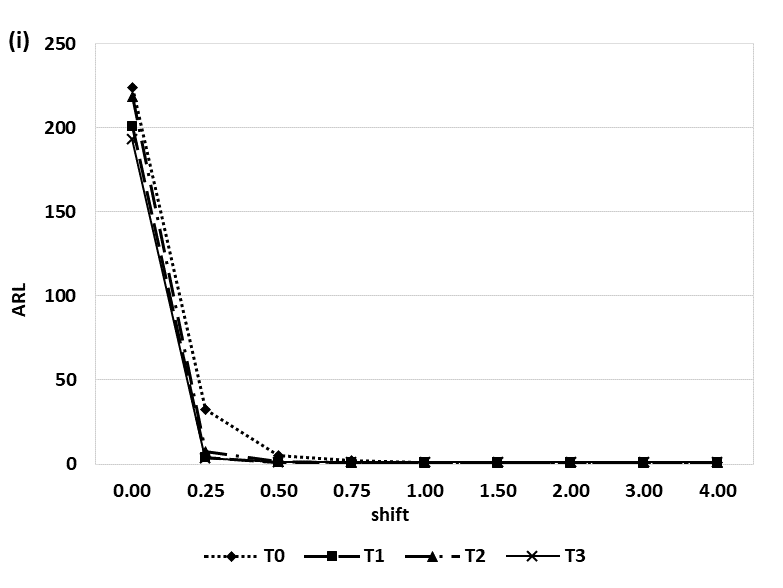}\hspace{0.1cm}&
		\includegraphics[width=8cm, height=6cm]{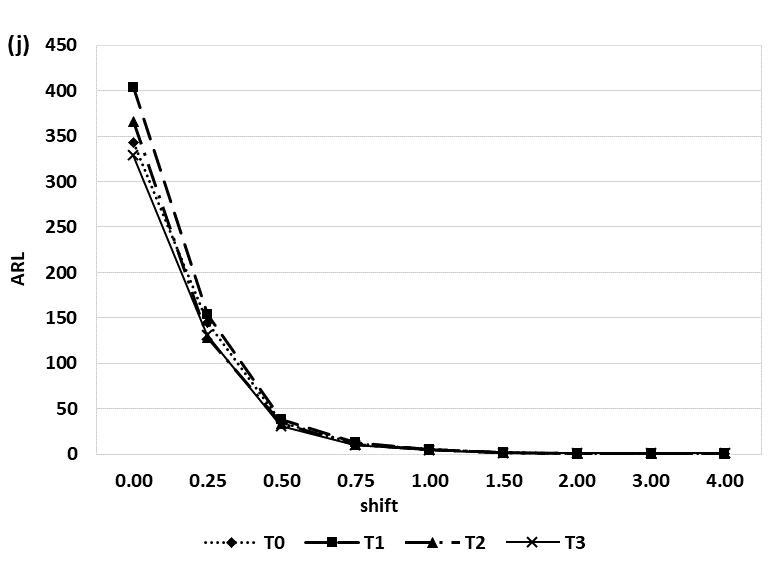}\hspace{0.1cm}\\
		
		\includegraphics[width=8cm, height=6cm]{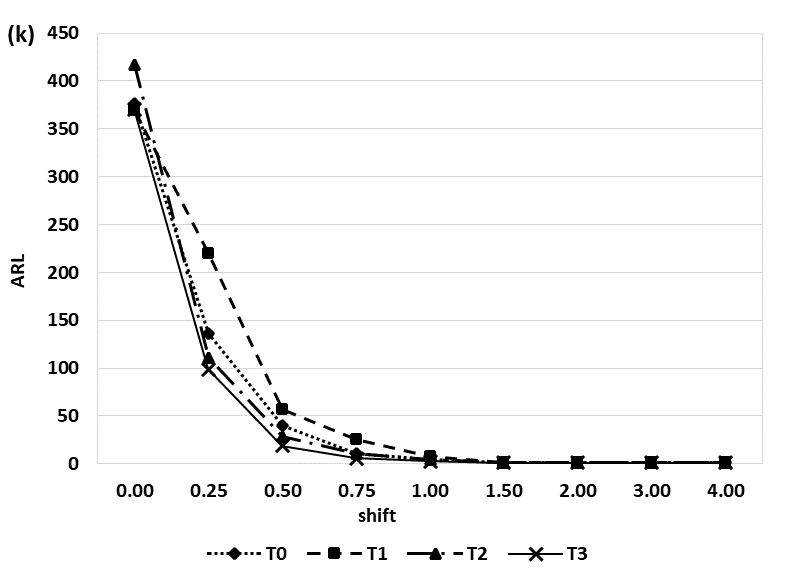}\hspace{0.1cm}&
		\includegraphics[width=8cm, height=6cm]{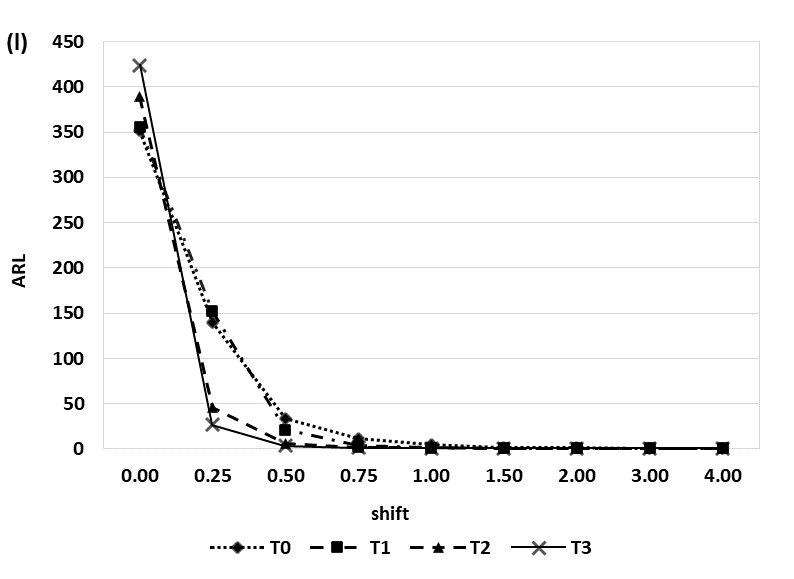}\hspace{0.1cm}\\
		
	\end{tabular}
	\caption{ARL curve of $T_0-T_3$ at: (g) ARL=200, n=15, $l_{xy}=0.30$; (h) ARL=200, n=15, $l_{xy}=0.60$; (i) ARL=200, n=15, $l_{xy}=0.90$; (j) ARL=371, n=5, $l_{xy}=0.30$; (k) ARL=371, n=5, $l_{xy}=0.60$; (l) ARL=371, n=5, $l_{xy}=0.90$.}
	\label{fig3}
\end{figure*} 

\begin{figure*}[h]
	\centering
	\begin{tabular}{c c}
		\includegraphics[width=8cm, height=6cm]{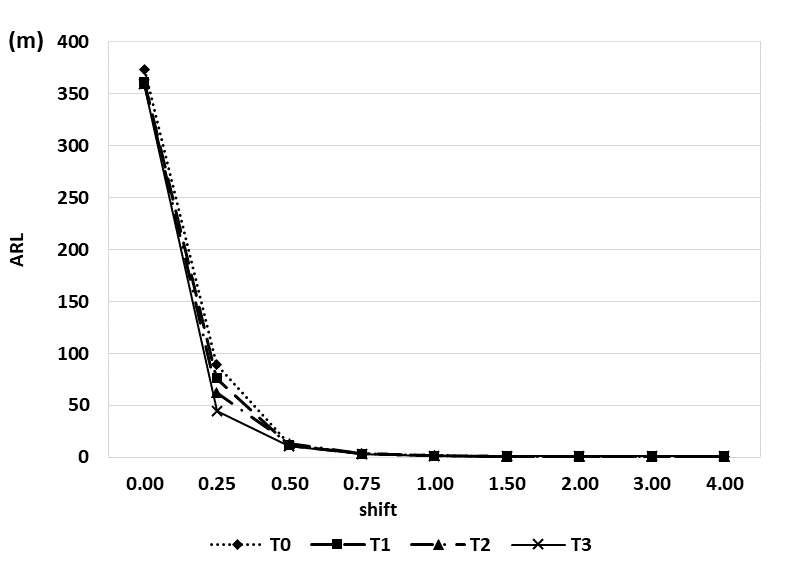}\hspace{0.1cm}&
		\includegraphics[width=8cm, height=6cm]{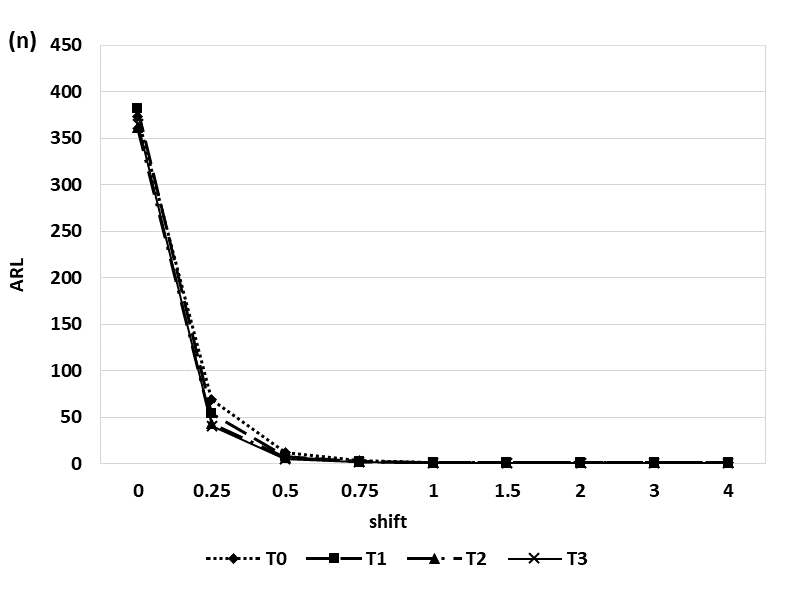}\hspace{0.1cm}\\
		
		\includegraphics[width=8cm, height=6cm]{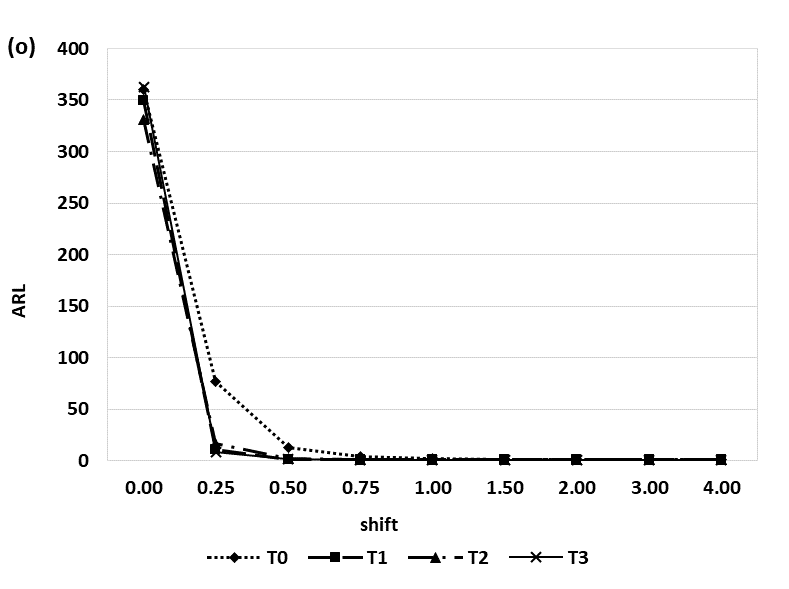}\hspace{0.1cm}&
		\includegraphics[width=8cm, height=6cm]{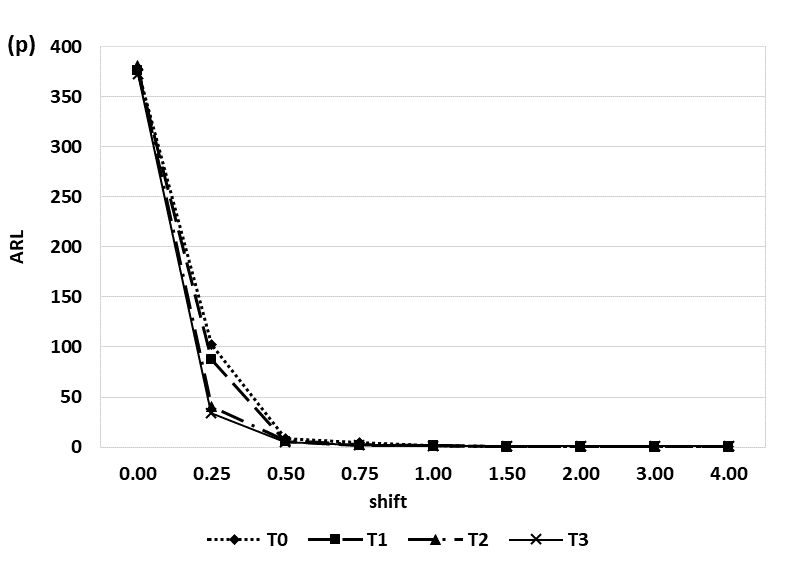}\hspace{0.1cm}\\
		
		\includegraphics[width=8cm, height=6cm]{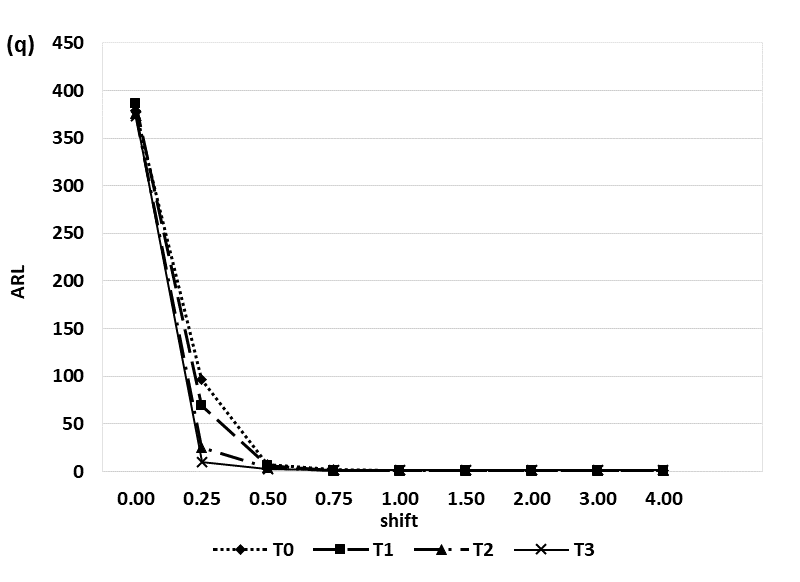}\hspace{0.1cm}&
		\includegraphics[width=8cm, height=6cm]{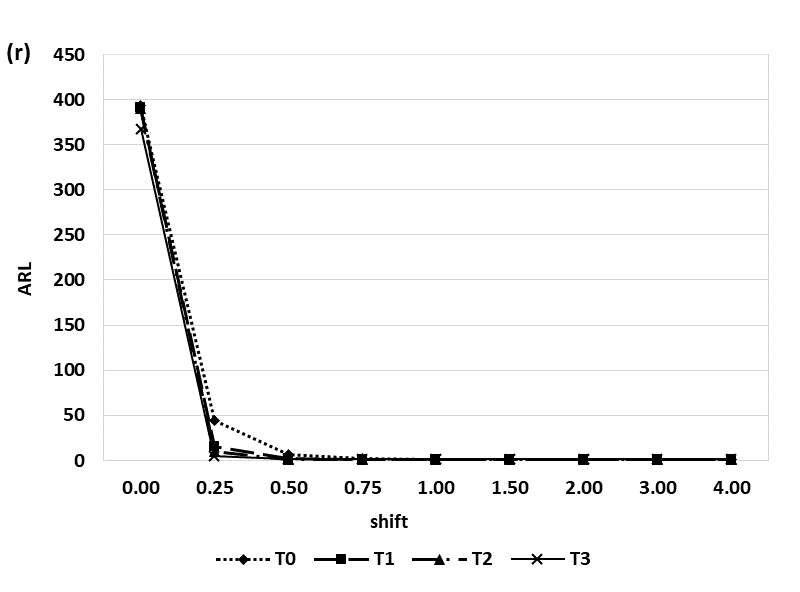}\hspace{0.1cm}\\
		
	\end{tabular}
	\caption{ARL curve of $T_0-T_3$ at: (m) ARL=371, n=10, $l_{xy}=0.30$; (n) ARL=371, n=10, $l_{xy}=0.60$; (o) ARL=371, n=10, $l_{xy}=0.90$; (p) ARL=371, n=15, $l_{xy}=0.30$; (q) ARL=371, n=15, $l_{xy}=0.60$; (r) ARL=371, n=15, $l_{xy}=0.90$.}
	\label{fig4}
\end{figure*}

\begin{figure*}[h]
	\centering
	\begin{tabular}{c c}
		\includegraphics[width=8cm, height=6cm]{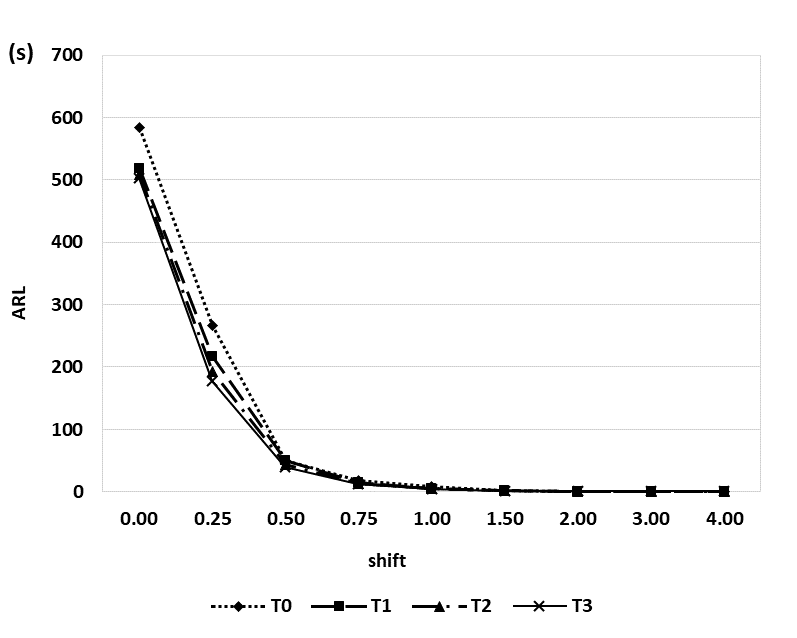}\hspace{0.1cm}&
		\includegraphics[width=8cm, height=6cm]{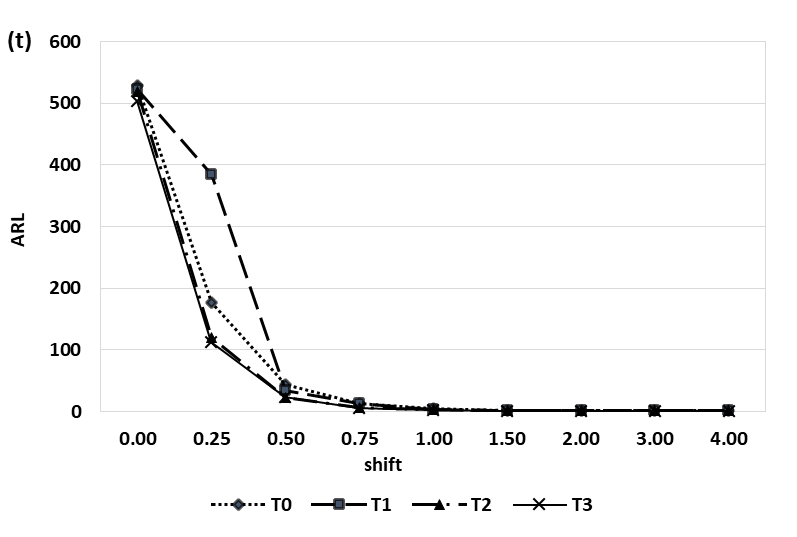}\hspace{0.1cm}\\
		
		\includegraphics[width=8cm, height=6cm]{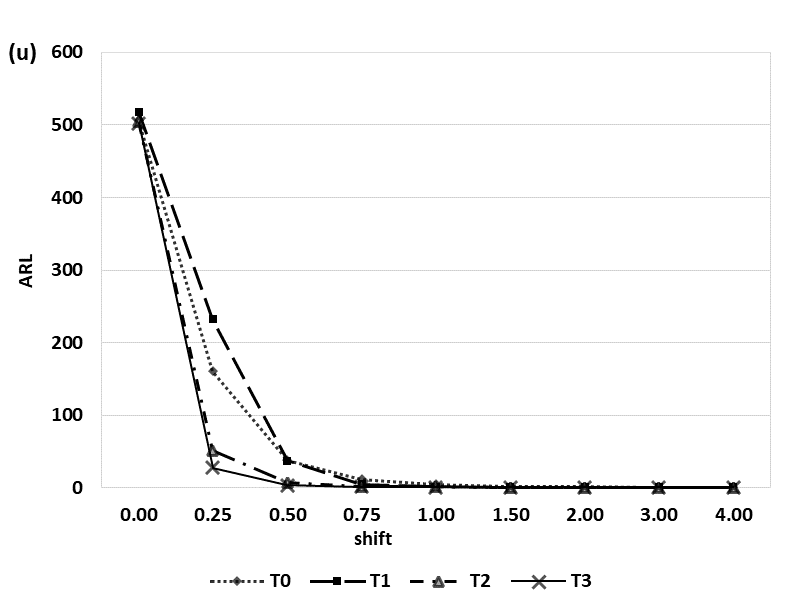}\hspace{0.1cm}&
		\includegraphics[width=8cm, height=6cm]{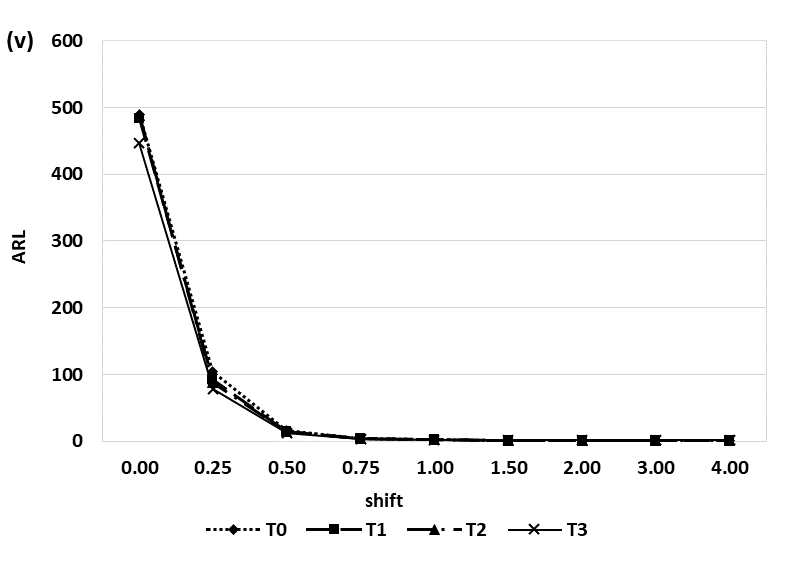}\hspace{0.1cm}\\
		
		\includegraphics[width=8cm, height=6cm]{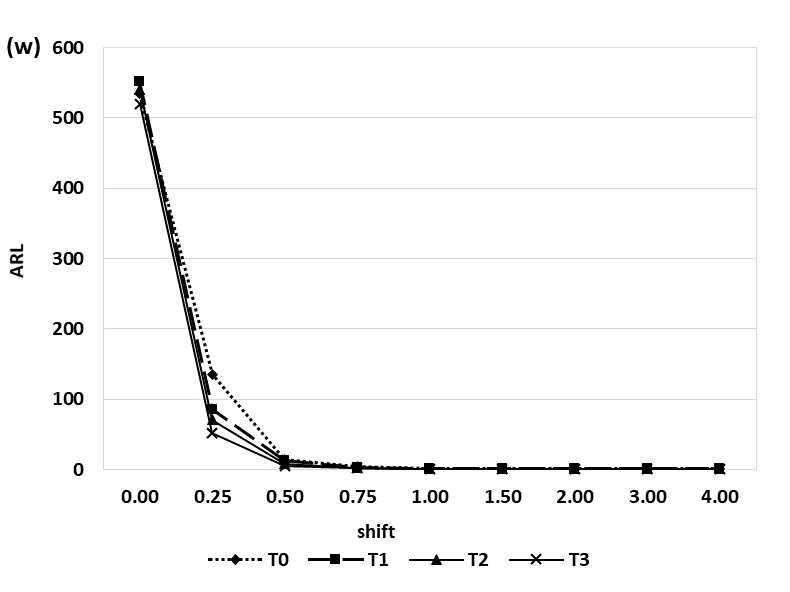}\hspace{0.1cm}&
		\includegraphics[width=8cm, height=6cm]{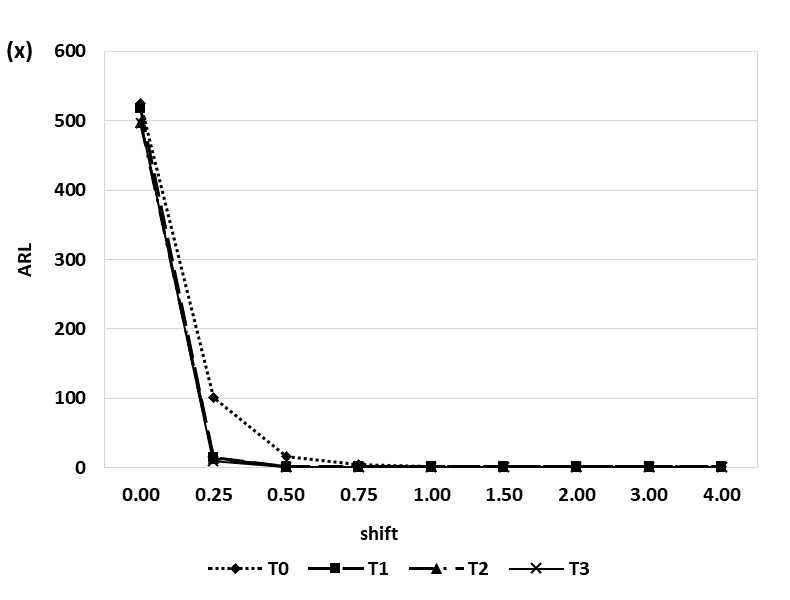}\hspace{0.1cm}\\
		
	\end{tabular}
	\caption{ARL curve of $T_0-T_3$ at: (s) ARL=500, n=5, $l_{xy}=0.30$; (t) ARL=500, n=5, $l_{xy}=0.60$; (u) ARL=500, n=5, $l_{xy}=0.90$; (v) ARL=500, n=10, $l_{xy}=0.30$; (w) ARL=500, n=10, $l_{xy}=0.60$; (x) ARL=500, n=10, $l_{xy}=0.90$.}
	\label{fig5}
\end{figure*}

\begin{figure*}[h]
	\centering
	\begin{tabular}{c c}
		\includegraphics[width=8cm, height=6cm]{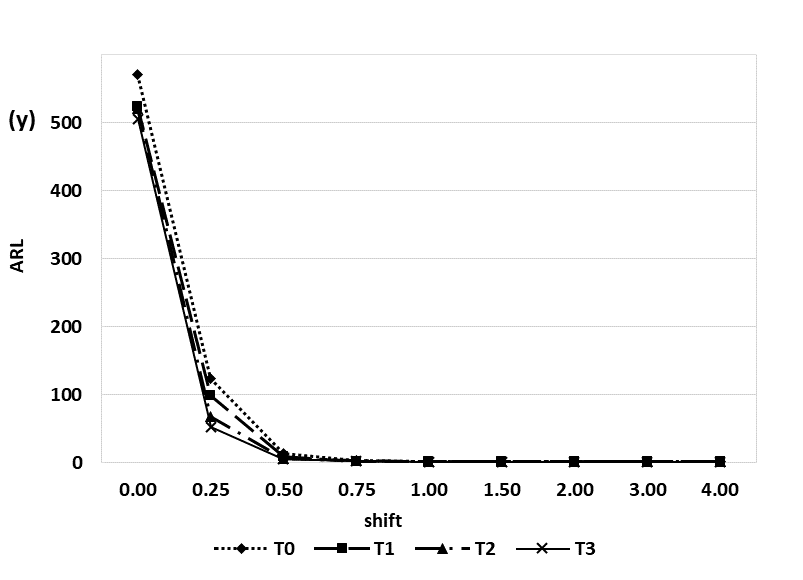}\hspace{0.1cm}&
		\includegraphics[width=8cm, height=6cm]{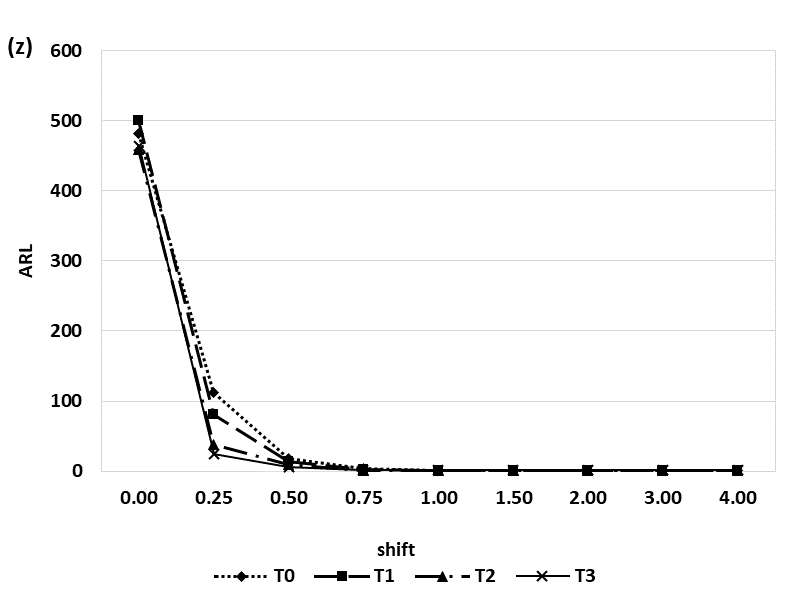}\hspace{0.1cm}\\
	\end{tabular}		
	\includegraphics[width=8cm, height=6cm]{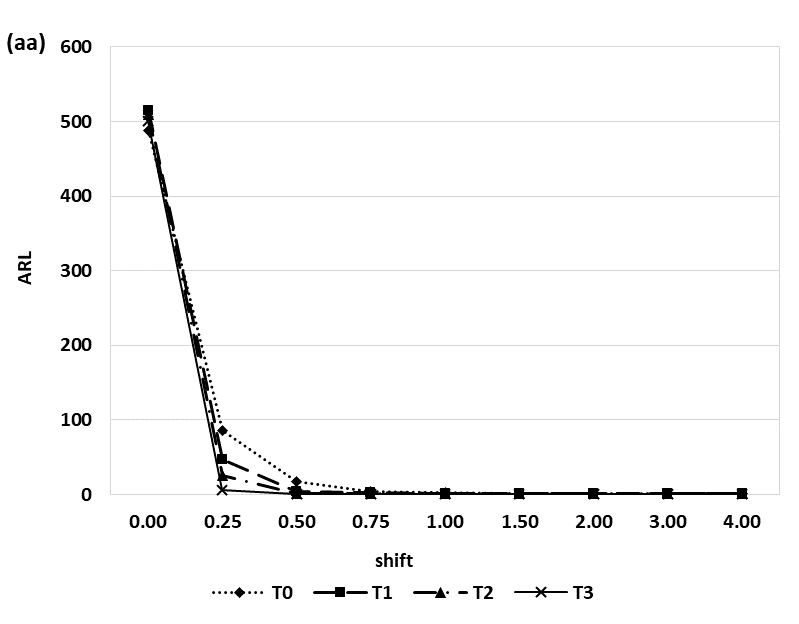}\hspace{0.1cm}
	
	\caption{ARL curve of $T_0-T_3$ at: (y) ARL=500, n=15, $l_{xy}=0.30$; (z) ARL=500, n=15, $l_{xy}=0.60$; (aa) ARL=500, n=15, $l_{xy}=0.90$.}
	\label{fig6}
\end{figure*}

\begin{table*}[h]
	\centering
	\caption{The values of $T_0, T_1, T_2, T_3$ for $\rho_{xy}=0.90$.}
	\begin{tabular}{lccccccccccccccc}\hline
		n=5&&Arl=500&&&&&Arl=371&&&&&Arl=200&&&\\
		\cline{1-1}\cline{3-6}\cline{8-11}\cline{13-16}
		Subgroup&&$T_0$&$T_1$&$T_2$&$T_3$&&$T_0$&$T_1$&$T_2$&$T_3$&&$T_0$&$T_1$&$T_2$&$T_3$\\\hline
		1&&4.623&5.165&6.425&5.052&&5.025&5.158&6.713&5.103&&5.420&5.228&7.008&5.189\\
		2&&4.424&5.294&6.461&5.271&&4.869&4.955&6.477&4.926&&5.216&4.759&6.613&4.853\\
		3&&4.442&5.058&6.208&4.925&&4.894&4.889&6.449&4.873&&5.085&5.415&6.972&5.375\\
		4&&5.016&4.948&6.585&4.953&&5.617&4.922&6.963&5.019&&4.408&4.788&6.161&4.894\\
		5&&4.993&4.789&6.455&4.812&&5.139&4.961&6.675&4.977&&5.425&4.891&6.726&4.856\\
		6&&4.595&5.081&6.388&5.028&&4.758&4.944&6.405&4.921&&4.733&4.926&6.347&4.866\\
		7&&4.266&5.373&6.288&5.190&&4.876&4.857&6.415&4.844&&4.199&4.763&5.964&4.823\\
		8&&5.236&4.905&6.667&4.904&&5.187&5.027&6.710&4.987&&5.412&5.198&6.962&5.140\\
		9&&5.065&5.036&6.658&5.008&&5.800&5.315&7.171&5.154&&4.874&4.748&6.318&4.728\\
		10&&4.791&4.510&6.144&4.573&&5.020&4.907&6.538&4.893&&5.357&4.719&6.476&4.605\\
		11&&5.924&5.108&7.269&5.197&&5.570&4.697&6.869&4.937&&4.958&4.932&6.525&4.920\\
		12&&4.458&4.786&6.178&4.873&&5.486&5.109&6.878&4.997&&5.129&4.638&6.427&4.687\\
		13&&4.590&4.923&6.271&4.884&&4.696&5.116&6.567&5.172&&4.651&4.995&6.361&4.949\\
		14&&5.436&5.210&6.862&5.009&&4.966&4.867&6.484&4.866&&5.723&4.574&6.927&4.918\\
		15&&5.107&5.151&6.792&5.140&&4.520&4.869&6.207&4.859&&4.753&5.277&6.634&5.211\\
		16&&4.290&4.904&6.048&4.851&&4.793&4.887&6.361&4.839&&4.332&4.880&6.198&5.010\\
		17&&5.513&5.095&6.759&4.843&&5.496&5.100&6.916&5.035&&4.587&4.640&6.070&4.632\\
		18&&4.908&5.089&6.649&5.109&&5.116&5.086&6.748&5.081&&4.516&4.729&6.144&4.782\\
		19&&5.407&5.088&6.925&5.100&&4.505&4.755&6.053&4.675&&5.200&5.259&6.922&5.232\\
		20&&4.423&4.902&6.266&5.018&&5.636&5.371&7.267&5.358&&4.926&5.040&6.647&5.093\\
		21&&4.613&4.875&6.237&4.823&&4.527&5.149&6.182&4.821&&3.924&4.680&5.732&4.779\\
		22&&5.113&5.178&6.835&5.188&&4.984&5.155&6.722&5.143&&4.778&4.595&6.134&4.569\\
		23&&4.735&4.928&6.387&4.914&&4.782&4.869&6.379&4.869&&4.602&5.046&6.374&5.004\\
		24&&5.470&5.255&6.963&5.106&&4.813&4.462&6.119&4.526&&5.666&5.092&7.152&5.207\\
		25&&5.402&4.868&6.848&5.013&&4.935&5.307&6.641&5.080&&5.066&5.299&6.866&5.259\\
		26&&5.363&4.857&6.816&5.000&&4.827&4.686&6.265&4.695&&5.563&5.104&7.097&5.205\\
		27&&4.950&5.159&6.719&5.165&&6.285&5.320&7.434&5.198&&5.415&4.891&6.711&4.845\\
		28&&5.672&4.753&6.878&4.890&&5.245&4.083&6.432&4.620&&5.055&5.377&6.893&5.299\\
		29&&4.905&5.187&6.727&5.208&&5.531&5.345&7.185&5.326&&4.933&5.054&6.558&4.979\\
		30&&4.278&4.396&5.709&4.419&&5.085&4.931&6.561&4.877&&5.627&5.095&6.946&4.993\\
		31&&6.596&6.018&8.287&5.976&&5.851&6.134&7.990&6.058&&6.302&6.296&8.409&6.264\\
		32&&6.213&6.253&8.366&6.267&&6.613&5.963&8.388&6.078&&5.648&5.700&7.557&5.686\\
		33&&5.175&6.111&7.605&6.081&&6.859&5.859&8.367&5.939&&5.549&5.965&7.845&6.091\\
		34&&6.014&6.256&8.250&6.255&&6.092&5.996&8.034&5.963&&5.876&5.696&7.636&5.638\\
		35&&6.131&6.280&8.342&6.289&&6.023&5.980&7.974&5.936&&5.449&5.605&7.447&5.687\\
		36&&5.873&6.277&8.217&6.307&&5.537&5.642&7.453&5.636&&6.412&6.261&8.374&6.166\\
		37&&5.280&6.090&7.809&6.249&&5.558&5.853&7.620&5.818&&5.606&6.237&7.755&5.945\\
		38&&5.559&5.752&7.651&5.854&&5.427&5.604&7.350&5.587&&6.671&6.290&8.320&5.975\\
		39&&5.774&5.791&7.699&5.771&&5.568&5.733&7.534&5.710&&6.825&6.573&8.764&6.379\\
		40&&6.221&5.784&7.881&5.722&&5.482&5.963&7.647&5.903&&6.347&5.974&8.089&5.886\\
		41&&6.175&6.054&8.147&6.043&&5.700&5.890&7.684&5.801&&7.021&6.452&8.690&6.210\\
		42&&6.146&5.867&8.004&5.899&&6.044&6.013&8.035&5.993&&5.988&5.708&7.762&5.716\\
		43&&5.720&5.668&7.566&5.651&&6.598&6.291&8.493&6.199&&6.203&6.080&8.188&6.073\\
		44&&6.341&6.238&8.383&6.214&&5.444&5.682&7.481&5.732&&5.920&6.106&8.096&6.137\\
		45&&5.759&5.871&7.774&5.868&&5.345&5.905&7.788&6.173&&5.859&6.017&7.969&6.029\\
		46&&6.172&6.135&8.183&6.085&&5.505&5.893&7.757&6.017&&6.410&6.261&8.441&6.240\\
		47&&6.496&6.216&8.375&6.122&&6.131&6.260&8.306&6.248&&5.922&5.961&7.933&5.949\\
		48&&6.453&5.997&8.148&5.896&&6.725&5.935&8.440&6.079&&6.339&5.908&8.134&5.939\\
		49&&6.568&6.148&8.360&6.069&&6.723&6.228&8.270&5.897&&5.042&5.543&7.292&5.801\\
		50&&5.320&5.811&7.506&5.849&&6.030&5.905&7.927&5.879&&4.553&5.936&7.140&6.044\\
		\hline
	\end{tabular}%
	\label{tab1}%
\end{table*}%

\begin{table*}[h]
	\centering
	\caption{The values of $T_0, T_1, T_2, T_3$ for $\rho_{xy}=0.60$ and $\rho_{xy}=0.90$.}
	\begin{tabular}{lcccccccccc}\hline
		&&n=10&$\rho_{xy}=0.60$&&&&n=15&$\rho_{xy}=0.90$&&\\
		\cline{3-6}\cline{8-11}
		Subgroup&&Arl=500&&&&&Arl=371&&&\\
		\cline{3-6}\cline{8-11}
		&&$T_0$&$T_1$&$T_2$&$T_3$&&$T_0$&$T_1$&$T_2$&$T_3$\\\hline
		1&&4.676&4.843&5.263&4.900&&4.915&6.397&4.919&4.919\\
		2&&4.737&4.965&5.319&4.948&&4.988&6.460&5.038&5.038\\
		3&&4.845&5.066&5.457&5.077&&5.152&7.014&5.099&5.099\\
		4&&4.551&4.534&4.889&4.516&&5.057&6.773&5.081&.081\\
		5&&5.473&5.166&5.678&5.201&&4.821&6.280&4.897&4.897\\
		6&&5.604&5.599&6.059&5.578&&4.663&6.060&4.651&4.651\\
		7&&4.882&4.971&5.399&5.004&&5.139&7.282&5.157&5.157\\
		8&&4.918&5.266&5.590&5.207&&5.125&6.824&5.345&5.345\\
		9&&5.320&5.184&5.608&5.150&&4.790&7.120&5.049&5.049\\
		10&&4.892&4.870&5.255&4.846&&5.307&6.820&5.286&5.286\\
		11&&4.727&5.091&5.462&5.112&&5.049&6.645&5.058&5.058\\
		12&&5.253&5.102&5.571&5.121&&4.678&6.096&4.660&4.66\\
		13&&4.437&4.572&4.885&4.536&&5.174&6.640&5.149&5.149\\
		14&&5.038&5.064&5.470&5.050&&5.380&7.106&5.363&5.363\\
		15&&5.372&5.123&5.514&5.045&&5.170&6.854&5.140&5.14\\
		16&&5.173&5.047&5.442&4.998&&4.917&6.411&4.988&4.988\\
		17&&5.057&5.338&5.649&5.240&&5.188&6.930&5.139&5.139\\
		18&&4.897&5.086&5.479&5.089&&5.206&6.742&5.244&5.244\\
		19&&5.056&5.102&5.546&5.129&&5.045&6.680&5.005&5.005\\
		20&&4.515&4.805&5.159&4.824&&5.150&7.021&5.154&5.154\\
		21&&4.582&4.765&5.131&4.775&&5.243&7.019&5.146&5.146\\
		22&&5.023&4.813&5.204&4.771&&5.319&6.754&5.309&5.309\\
		23&&4.876&4.832&5.222&4.814&&5.238&6.936&5.217&5.217\\
		24&&5.219&5.023&5.536&5.089&&4.570&5.982&4.554&4.554\\
		25&&4.838&4.724&5.111&4.703&&4.887&5.930&4.733&4.733\\
		26&&5.238&5.180&5.623&5.178&&4.845&6.170&4.746&4.746\\
		27&&4.770&4.790&5.163&4.771&&5.234&6.474&5.177&5.177\\
		28&&4.782&5.065&5.384&5.010&&4.824&6.260&4.776&4.776\\
		29&&5.035&5.154&5.531&5.116&&5.071&6.454&5.025&5.025\\
		30&&5.122&5.231&5.628&5.205&&4.961&6.196&4.884&4.884\\
		31&&5.986&5.948&6.438&5.916&&5.784&7.482&5.748&5.748\\
		32&&5.867&5.794&6.242&5.730&&5.934&7.736&5.997&5.997\\
		33&&5.817&6.028&6.508&6.018&&5.659&7.537&5.867&5.867\\
		34&&6.015&5.928&6.456&5.930&&5.997&8.061&5.996&5.996\\
		35&&6.070&6.004&6.484&5.952&&5.577&7.278&5.826&5.826\\
		36&&6.197&6.324&6.914&6.382&&6.219&8.506&6.077&6.077\\
		37&&6.178&5.926&6.386&5.839&&5.605&7.883&5.618&5.618\\
		38&&5.874&6.268&6.774&6.296&&6.171&8.139&5.926&5.926\\
		39&&5.734&5.869&6.407&5.926&&6.472&8.257&6.283&6.283\\
		40&&5.958&5.864&6.310&5.788&&6.038&8.126&5.901&5.901\\
		41&&6.389&6.569&7.195&6.647&&5.866&8.045&5.776&5.776\\
		42&&6.254&6.184&6.704&6.153&&6.079&7.875&6.093&6.093\\
		43&&6.468&6.311&6.813&6.238&&5.792&7.753&5.788&5.788\\
		44&&6.046&6.368&7.002&6.508&&6.063&8.100&6.036&6.036\\
		45&&5.734&5.654&6.125&5.627&&5.872&7.703&6.077&6.077\\
		46&&5.995&6.124&6.588&6.072&&5.636&8.272&5.965&5.965\\
		47&&6.127&5.900&6.361&5.819&&6.228&8.165&6.274&6.274\\
		48&&5.639&5.782&6.337&5.869&&5.591&8.461&6.060&6.060\\
		49&&6.325&6.266&6.748&6.189&&5.812&7.867&5.718&5.718\\
		50&&6.728&6.525&6.908&6.307&&6.120&7.384&6.195&6.195\\
		\hline
	\end{tabular}%
	\label{tab2}%
\end{table*}%

\end{document}